\documentclass[12pt]{amsart}
\usepackage{amssymb}
\usepackage[all]{xy}

\long\def\forget#1\forgotten{}
\newcommand{\issuenumber}{16}
\newcommand{\issuemonth}{March}
\newcommand{\issueyear}{2006}

\newtheorem{thm}{Theorem}[section]
\newtheorem{prob}[thm]{Problem}

\newtheorem{issue}{Issue}

\theoremstyle{definition}

\theoremstyle{remark}

\newcommand{\ed}{
\general\end{document}}

\newcommand{\fb}{\mathfrak{b}}

\newcommand{\fd}{\mathfrak{d}}
\newcommand{\fp}{\mathfrak{p}}

\newcommand{\NON}{{\mathsf   {NON}}}
\newcommand{\COF}{{\mathsf   {COF}}}

\newcommand{\M}{\mathcal{M}}

\newcommand{\cov}{\mathsf{cov}}

\newcommand{\CH}{the Continuum Hypothesis}
\newcommand{\R}{\mathbb{R}}

\newcommand{\fo}{\mathfrak{od}}

\renewcommand{\b}{\mathfrak{b}}

\renewcommand{\split}{\mathsf{Split}}
\newcommand{\bq}{\begin{quote}}
\newcommand{\eq}{\end{quote}}
\renewcommand{\O}{\mathcal{O}}
\newcommand{\B}{\mathcal{B}}
\newcommand{\BG}{\B_\Gamma}

\newcommand{\BO}{\B_\Omega}

\newcommand{\sone}{\mathsf{S}_1}    \newcommand{\sfin}{\mathsf{S}_{fin}}

\newcommand{\ufin}{\mathsf{U}_{fin}}

\newcommand{\nin}{\not\in}

\newcommand{\sbst}{\subseteq}
\newcommand{\by}[2]{\par\hfill\emph{#1}, #2}
\newcommand{\nby}[1]{\par\hfill\emph{#1}}
\newcommand{\Tau}{\mathrm{T}}
\newcommand{\CE}{\textsc{CE}}

\newcommand{\be}{\begin{enumerate}}
\newcommand{\ee}{\end{enumerate}}
\newcommand{\bi}{\begin{itemize}}
\newcommand{\ei}{\end{itemize}}
\renewcommand{\i}{\item}

\newcommand{\general}{\small\vfill\par\noindent\hrulefill\par
\noindent\textbf{Previous issues.} The first issues of this
bulletin, which contain general information (first issue), basic
definitions, research announcements, and open problems (all
issues) are available online, on \arx{math.GN/$x$}, where $x$ is
\texttt{0301011}, \texttt{0302062}, \texttt{0303057},
\texttt{0304087}, \texttt{0305367}, \texttt{0312140},
\texttt{0401155}, \texttt{0403369}, \texttt{0406411},
\texttt{0409072}, \texttt{0412305}, \texttt{0503631},
\texttt{0508563}, \texttt{0509432}, and \texttt{0512275},
respectively, for issues number $1$ to $15$.\\[0.1cm]
\textbf{Contributions.}
Please submit your contributions (announcements, discussions, and open problems)
by e-mailing us. It is preferred to write them
in \LaTeX{}.
The authors are urged to use as standard notation as possible, or otherwise give
the definitions or a reference to where the notation is explained.
Contributions to this bulletin would not require any transfer of copyright,
and material presented here can be published elsewhere.\\[0.1cm]
\textbf{Subscription.}
To receive this bulletin (free) to your
e-mailbox, e-mail us:\\
{boaz.tsaban@weizmann.ac.il}
}

\newcommand{\nArxPaper}[5]{\subsection{#2}{#4}\par\hfill{\arx{#1}}\par\hfill\emph{#3}}

\newcommand{\nAMSPaper}[5]{\subsection{#2}{#4}\par\hfill{\texttt{#1}}\par\hfill\emph{#3}}

\newcommand{\arx}[1]{\texttt{http://arxiv.org/abs/#1}}
\newcommand{\url}[1]{\bq\texttt{#1}\eq}
\newcommand{\online}[1]{The paper is available online at \url{#1}}

\title[$\mathcal{SPM}$ Bulletin \textbf{\issuenumber} (\issuemonth{} \issueyear)]{%
$\mathcal{SPM}$ Bulletin\\[0.5cm]
Issue number \issuenumber: \issuemonth{} \issueyear{} \CE{}}

\begin{document}
\maketitle

\tableofcontents

\section{Editor's note}

Some good news in short:
\be
\i The published version of Moore's paper \emph{A solution to the $L$ space problem}
is now available online at
\url{http://www.ams.org/journal-getitem?pii=S0894-0347-05-00517-5}
\i Maharam's Problem was solved by Talagrand (see below).
\i In the last Jerusalem Logic Seminar, Shelah reported on a new try to solve the minimal tower problem (see Issue 5).
\i The BEST meeting and the satellite meeting on SPM are very nearby.
\ee

\medskip

Contributions to the next issue are, as always, welcome.

\medskip

\by{Boaz Tsaban}{boaz.tsaban@weizmann.ac.il}

\hfill \texttt{http://www.cs.biu.ac.il/\~{}tsaban}

\section{Research announcements}

\nAMSPaper{http://www.ams.org/journal-getitem?pii=S1079-6762-05-00152-6}
{Hurewicz-like tests for Borel subsets of the plane}
{Dominique Lecomte}
{Let $\xi \geq 1$ be a countable ordinal. We study the Borel
subsets of the plane that can be made $\boldsymbol \Pi ^{0}_{\xi }$
by refining the Polish topology
on the real line. These sets are called
potentially $\boldsymbol \Pi ^{0}_{\xi }$. We give a
Hurewicz-like test to recognize potentially $\boldsymbol \Pi ^{0}_{\xi }$ sets.
}

\nArxPaper{math.GN/0512553}
{Ordered Spaces, Metric Preimages, and Function Algebras}
{Kenneth Kunen}
{We prove some results about compact Hausdorff spaces which have scattered-to-one maps onto compact metric spaces,
along with two types of consequences of these results for compact LOTSes (totally ordered spaces).
The first type shows that many products of $n+1$ LOTSes cannot be embedded into any product of $n$ LOTSes.
The second involves the Complex Stone-Weierstrass Property (CSWP), which is the complex version of the
Stone-Weierstrass Theorem. If X is a compact subspace of a product of three LOTSes, then $X$ has the CSWP if and only if
$X$ has no subspace homeomorphic to the Cantor set.
}

\nArxPaper{math.LO/0512546}
{On the independence of a generalized statement of Egoroff's theorem from ZFC, after T.\ Weiss}
{Roberto Pinciroli}
{We consider a generalized version (GES) of the wellknown Severini-Egoroff theorem in real analysis,
first shown to be undecidable in ZFC by Tomasz Weiss.
This independence is easily derived from suitable hypotheses on some cardinal characteristics of the continuum
like $\fb$ and $\fo$, the latter being the least cardinality of a subset of $[0,1]$ having full outer measure.}

\nArxPaper{math.GN/0512564}
{Forty annotated questions about large topological groups}
{Vladimir Pestov}
{This is a selection of open problems dealing with ``large'' (non locally
compact) topological groups and concerning extreme amenability
(fixed point on compacta property), oscillation stability,
universal minimal flows and other aspects of universality, and
unitary representations.}

\nAMSPaper{http://www.ams.org/journal-getitem?pii=S0002-9939-05-08206-7}
{Strong compactness and a partition property}
{Pierre Matet}
{We show that if $\operatorname{Part}(\kappa,\lambda)$ holds for
every $\lambda\ge\kappa$,
then $\kappa$ is strongly compact.
}

\nAMSPaper{math.LO/0512626}
{Countable Borel equivalence relations and quotient Borel spaces}
{Roberto Pinciroli}
{We consider countable Borel equivalence relations on quotient Borel spaces.
We prove a generalization of the Feldman-Moore representation theorem, but provide some examples
showing that other very simple properties of countable equivalence relations on standard Borel spaces
may fail in the context of nonsmooth quotients.}

\nArxPaper{math.LO/0601083}
{Decisive creatures and large continuum}
{Jakob Kellner and Saharon Shelah}
{
\newcommand{\ho}{\ensuremath{^{\omega}}}                  
\newcommand{\mycfa}{c^{\forall}}
\newcommand{\al}[1]{\ensuremath{{\aleph_{#1}}} }          
\newcommand{\myc}{c^{\exists}}
  For $f,g\in\omega\ho$ let $\mycfa_{f,g}$ be the minimal
  number of uniform trees with $g$-splitting needed to
  $\forall^\infty$-cover a uniform tree with $f$-splitting.
  $\myc_{f,g}$ is the dual notion for the $\exists^\infty$-cover.

  Assuming CH and given
  $\al1$ many (sufficiently different) pairs $(f_\epsilon,g_\epsilon)$ and
  cardinals $\kappa_\epsilon$ such that $\kappa_\epsilon^\al0=\kappa_\epsilon$,
  we construct a partial order forcing that
  $\myc_{f_\epsilon,g_\epsilon}=\mycfa_{f_\epsilon,g_\epsilon}=\kappa_\epsilon$.

  For this, we introduce a countable support semiproduct of
  decisive creatures with bigness and halving. This semiproduct
  satisfies fusion, pure decision and continuous reading of names.
}

\nArxPaper{math.LO/0601087}
{Models of real-valued measurability}
{Saka\'{e} Fuchino, Noam Greenberg, and Saharon Shelah}
{Solovay's random-real forcing (1971) is the standard way of producing
real-valued measurable cardinals. Following questions of Fremlin, by giving a new
construction, we show that there are combinatorial, measure-theoretic properties of
Solovay's model that do not follow from the existence of real-valued measurability.}

\nAMSPaper{http://www.ams.org/journal-getitem?pii=S0002-9939-06-08433-4}
{Hausdorff ultrafilters}
{Mauro Di Nasso and Marco Forti}
{We give the name \emph{Hausdorff} to those ultrafilters that provide
ultrapowers whose natural topology ($S$-topology) is Hausdorff, e.g. selective
ultrafilters are Hausdorff. Here we give  necessary and sufficient conditions
for  product ultrafilters  to be Hausdorff. Moreover we show  that no regular
ultrafilter over the ``small'' uncountable cardinal $\mathfrak{u}$ can be Hausdorff.
($\mathfrak{u}$ is the least size of an ultrafilter basis on $\omega$.)  We focus on
countably incomplete ultrafilters, but our main results also hold for
$\kappa$-complete ultrafilters.
}

\nAMSPaper{http://www.ams.org/journal-getitem?pii=S0002-9947-06-03864-5}
{Block combinatorics}
{V.\ Farmaki and S.\ Negrepontis}
{In this paper we extend the block combinatorics partition theorems of
Hindman and Milliken-Taylor in the setting of the recursive system of the block
Schreier families $(\mathcal{B}^\xi)$, consisting of families defined for every
countable ordinal $\xi$.
Results contain
(a)~a block partition Ramsey theorem for every countable
ordinal $\xi$ (Hindman's Theorem corresponding to $\xi=1$, and the Milliken-Taylor
Theorem to $\xi$ a finite ordinal),
(b)~a countable ordinal form of the block Nash-Williams partition theorem,
and (c)~a countable ordinal block partition theorem for sets closed in the
infinite block analogue of Ellentuck's topology.}

\nArxPaper{math.FA/0601689}
{Maharam's problem}
{Michel Talagrand}
{We construct an exhaustive submeasure that is not equivalent to a measure.
This solves problems of J.\ von Neumann (1937) and D.\ Maharam (1947).}

\nAMSPaper{http://www.ams.org/journal-getitem?pii=S0002-9939-06-08189-5}
{Universality of uniform Eberlein compacta}
{Mirna Dzamonja}
{We prove that if $\mu^+
<\lambda={\rm cf}(\lambda)<\mu^{\aleph_0}$ for some regular
$\mu>2^{\aleph_0}$, then there is no family of
less than $\mu^{\aleph_0}$ c-algebras of size $\lambda$ which are
jointly universal for c-algebras of size $\lambda$. On the
other hand, it is consistent to have a cardinal $\lambda\ge \aleph_1$
as large as desired and satisfying
$\lambda^{<\lambda}=\lambda$ and $2^{\lambda^+}>\lambda^{++}$, while there
are $\lambda^{++}$ c-algebras of size $\lambda^+$ that are jointly
universal for c-algebras of size $\lambda^+$. Consequently, by the
known results of M. Bell, it is consistent that there is $\lambda$ as
in the last statement and $\lambda^{++}$ uniform Eberlein compacta
of weight $\lambda^+$ such that at least one among them maps onto any Eberlein compact
of weight $\lambda^+$ (we call such a family universal).
The only positive universality results
for Eberlein compacta known previously required the relevant instance of $GCH$
to hold. These results complete the answer to a question of Y. Benyamini,
M. E. Rudin and M. Wage from 1977 who asked if there always was a universal
uniform Eberlein compact of a given weight.}

\nArxPaper{math.FA/0602628}
{Linearly ordered compacta and Banach spaces with a projectional resolution of the identity}
{Wieslaw Kubis}
{We construct a compact linearly ordered space $K$ of weight aleph one, such
that the space $C(K)$ is not isomorphic to a Banach space with a projectional
resolution of the identity, while on the other hand, $K$ is a continuous image
of a Valdivia compact and every separable subspace of $C(K)$ is contained in a
1-complemented separable subspace. This answers two questions due to O.\ Kalenda
and V.\ Montesinos.}

\subsection{Steinhaus Sets and Jackson Sets}
We prove that there does not exist a subset of the plane $S$ that meets every isometric copy of the
vertices of the unit square in exactly one point. We give a complete characterization of all three
point subsets $F$ of the reals such that there does not exists a set of reals $S$ which meets every
isometric copy of $F$ in exactly one point. A finite set $X$ in the plane is Jackson iff for every subset
$S$ of the plane there exists an isometric copy $Y$ of $X$ such that $Y$ does not meets $S$ in exactly
one point. These results are related to the open problem (Steve Jackson):
Is every finite set $X$ in the plane of two or more points Jackson?

\url{http://www.math.wisc.edu/~miller/res/jack.pdf}

\nby{Su Gao, Arnold W.\ Miller, and William A.\ R.\ Weiss}

\section{Problem of the Issue}

For a set of reals $X$, $\BO$ denotes the collection of all countable Borel $\omega$-covers of $X$,
and $\BO$ denotes the collection of all countable Borel $\gamma$-covers of $X$.
$\sone(\BO,\BG)$ is equivalent to $\binom{\BO}{\BG}$, that is,
``every countable Borel $\omega$-cover of $X$ contains a $\gamma$-cover
of $X$''.

The following problem is posed
in Miller's paper \emph{A Nonhereditary Borel-cover $\gamma$-set}
(Real Analysis Exchange \textbf{29} (2003/4), 601--606).

\begin{prob}
Does Martin's Axiom imply the existence of an uncountable set of reals
satisfying $\sone(\BO,\BG)$?
\end{prob}

It is known that \CH{} can be used for that (see e.g.\ Miller's mentioned
paper). There are unpublished notes of Todorcevic on this problem, contact
me directly for a copy.

\nby{Boaz Tsaban}

\newpage

\section{Problems from earlier issues}

\begin{issue}
Is $\binom{\Omega}{\Gamma}=\binom{\Omega}{\Tau}$?
\end{issue}

\begin{issue}
Is $\ufin(\Gamma,\Omega)=\sfin(\Gamma,\Omega)$?
And if not, does $\ufin(\Gamma,\Gamma)$ imply
$\sfin(\Gamma,\Omega)$?
\end{issue}

\stepcounter{issue}

\begin{issue}
Does $\sone(\Omega,\Tau)$ imply $\ufin(\Gamma,\Gamma)$?
\end{issue}

\begin{issue}
Is $\fp=\fp^*$? (See the definition of $\fp^*$ in that issue.)
\end{issue}

\begin{issue}
Does there exist (in ZFC) an uncountable set satisfying $\sone(\BG,\B)$?
\end{issue}

\stepcounter{issue}

\begin{issue}
Does $X \nin \NON(\M)$ and $Y\nin\mathsf{D}$ imply that
$X\cup Y\nin \COF(\M)$?
\end{issue}

\begin{issue}
Assume CH. Is $\split(\Lambda,\Lambda)$ preserved under finite unions?
\end{issue}

\begin{issue}
Is $\cov(\M)=\fo$? (See the definition of $\fo$ in that issue.)
\end{issue}

\begin{issue}
Does $\sone(\Gamma,\Gamma)$ always contain an element of cardinality $\b$?
\end{issue}

\begin{issue}
Could there be a Baire metric space $M$ of weight $\aleph_1$ and a partition
$\mathcal{U}$ of $M$ into $\aleph_1$ meager sets where for each ${\mathcal U}'\subset\mathcal U$,
$\bigcup {\mathcal U}'$ has the Baire property in $M$?
\end{issue}

\stepcounter{issue} 

\begin{issue}
Does there exist (in ZFC) a set of reals $X$ of cardinality $\fd$ such that all
finite powers of $X$ have Menger's property $\ufin(\O,\O)$?
\end{issue}

\begin{issue}
Can a Borel non-$\sigma$-compact group be generated by a Hurewicz subspace?
\end{issue}

\begin{issue}
Does MA imply the existence of an uncountable $X\sbst\R$
satisfying $\sone(\BO,\BG)$?
\end{issue}

\ed